\documentclass[11pt]{article}
\usepackage{graphicx}
\usepackage{amssymb,amsfonts,amsthm}
\usepackage{amsmath,amsthm,amssymb}
\usepackage{amsfonts}

\newtheorem{theorem}{Theorem}[section]
\newtheorem{lemma}[theorem]{Lemma}

\newtheorem{cy}[theorem]{Corollary}

\theoremstyle{definition}

\newtheorem{rk}[theorem]{Remark}

\newcounter{ppp}

\newcommand{\Ker}{\hbox{Ker}}

\begin{document}
\title{On identities in the products of group varieties}
\author{ N. S. Boatman, A. Yu. Olshanskii \thanks{The
author was supported in part by the NSF grant DMS 1161294}}
\maketitle

\begin{abstract}
Let ${\cal B}_n$ be the variety of groups satisfying the law $x^n=1$. It is proved that
for every sufficiently large prime $p$, say $p>10^{10}$, the product ${\cal B}_p{\cal B}_p$ cannot be
defined by a finite set of identities. This solves the problem formulated by
C.K. Gupta and A.N. Krasilnikov in 2003. We also find the axiomatic and the basis ranks of
the variety ${\cal B}_p{\cal B}_p$. For this goal,  we improve the estimate for the basis rank
of the product of group varieties obtained  by G. Baumslag, B.H. Neumann, H. Neumann and P.M. Neumann long ago.
\end{abstract}



\section{Introduction}

A variety of groups $\cal V$ is a class of groups given by a set of identical relations.
For example, the variety of abelian groups $\cal A$ is given by the identity $[x,y]=1$, and
the Burnside variety ${\cal B}_n$ contains all the groups with the identity
$x^n=1$. In general, a group identity looks like $v(x_1,\dots,x_n)=1$, where
$v(x_1,\dots,x_n)$ is a word, i.e., an element of a free group with a basis $x_1,\dots,x_n,\dots$.

For an arbitrary group $G$, the {\it verbal subgroup} $V(G)$ is generated by all the values $v(g_1,\dots, g_n)$ ($g_1, \dots, g_n\in G$) of
the left-hand sides of the identities $v(x_1,\dots,x_n)=1$ of the variety $\cal V$, and so $V(G)$ is the minimal normal
subgroup $N$ of $G$ such that $G/N\in \cal V$. For example, $V(G)$ is the derived subgroup of $G$
if $\cal V= A$ and $V(G)=G^n$ is generated by all the $n$-th powers of the elements from $G$ if $ {\cal V = B}_n.$

Given an arbitrary cardinal $m$, every variety $\cal V$ has a $\cal V$-free group $F_m(\cal V)$ of free rank $m$ that is isomorphic to
the factor group $F_m/V(F_m)$, where $F_m$ is the absolutely free group of rank $m$.

There are many  varieties which are
finitely based, that is,
can be defined by finite sets of identities.  For example, all nilpotent varieties \cite{L},
all metabelian and nilpotent-by-abelian varieties \cite{C, Kr}, and every variety generated by a single finite group \cite{OP} can be defined by a finite set of identities. We refer the reader to the book \cite{HN} and also to the surveys  \cite{BO} and \cite{GK}.
In the paper \cite{N} (and also in his Ph.D. thesis in 1935), B.Neumann asked
whether every group variety was finitely based. In 1969,
the negative answers were, in turn, given in the papers \cite{O70}, \cite{A} and \cite{VL}
based on different methods.

 A simplest example of an infinite system of group identities
which is not equivalent to a finite one was presented in \cite{B} and \cite{K73}:
\begin{equation}\label{B4B2}
x_1^8=1, (x_1^2x_2^2)^4=1,\dots, (x_1^2x_2^2\dots x_m^2)^4=1,\dots
\end{equation}
Note that a group $G$ satisfies the system of identities
(\ref{B4B2}) if and only if it has a normal subgroup $N$ satisfying the identity
$x^4=1$ such that the quotient $G/N$ satisfies the identity $x^2=1$. In other
words, the set (\ref{B4B2}) defines the product ${\cal B}_4{\cal B}_2$ of two
Burnside varieties ${\cal B}_4$ and ${\cal B}_2$. In general, the product $\cal UV$ of
group varieties $\cal U$ and $\cal V$ consists of all groups $G$ having a normal
subgroup $N\in \cal U$ such that $G/N\in \cal V$.

On one hand, an easy exercise shows that if positive integers $m$ and $n$ are relatively prime, then all the identities of the product ${\cal B}_m{\cal B}_n$ follow from the single identity $(x_1^nx_2^n)^m=1$. On the other hand, the
following generalization of the example (\ref{B4B2}) is obtained in \cite{K74}:
For any prime $p$, the product ${\cal B}_{n}{\cal C}$ is not finitely based if $n$ is divisible by $p^2$ and the variety $\cal C$ contains a cyclic group of order $p$, where $\cal C$
is not the variety of all groups. The hypothesis $p^2\mid n$ was used in an essential way in \cite{K74};
in particular, the case ${\cal B}_p{\cal B}_p$ remained unknown.
The authors of the survey \cite{GK} noticed that the variety ${\cal B}_p{\cal B}_p$
is finitely based for $p=2,3$, and they explicitly formulated
the question (Problem 2 in \cite{GK}): ``Is it true that for a large prime $p$ the variety
defined by the identities

\begin{equation}\label{BpBp}
x_1^{p^2}=1, (x_1^px_2^p)^p=1,\dots, (x_1^px_2^p\dots x_m^p)^p=1,\dots
\end{equation}
is not finitely based ?''

Here we answer this question.

\begin{theorem}\label{main} For every sufficiently large prime $p$, the
variety ${\cal B}_p{\cal B}_p$ is not finitely based, i.e. the system
(\ref{BpBp}) is not equivalent to a finite system of identities.
\end{theorem}

This theorem can be formulated in a stronger form. Recall that a group variety $\cal V$
has infinite {\it axiomatic  rank} if it cannot be given by any set of identities
(finite or infinite) depending on a finite set of variables. A variety $\cal V$ has
{\it basis rank} $d$ if it is generated by its $d$-generated groups, but is not generated
by $(d-1)$-generated groups, i.e. the minimal variety containing all $(d-1)$-generated
groups from $\cal V$ is strictly less than $\cal V$.

\begin{theorem}\label{rank} For every sufficiently large prime $p$,

(1) the variety ${\cal B}_p{\cal B}_p$ has infinite axiomatic rank;

(2) its basis rank equal to $2$.
\end{theorem}

The proof of statement (2) is based on Theorem \ref{UV} giving a new estimate of the basis ranks of
the products $\cal UV$ for various  $\cal U$ and $\cal V$.

\begin{rk} The methods of \cite{B, K73, K74} work for the varieties ${\cal B}_{p^2}{\cal B}_p$ (and many other
products of varieties) but cannot be applied to ${\cal B}_{p}{\cal B}_p$ for the following reason.
The authors of \cite{B, K73, K74} check identities in solvable groups of bounded length and prove, in fact, that
even the intersections ${\cal B}_{p^2}{\cal B}_p\cap \cal S$ have no finite bases of identities,
where $\cal S$ is a finitely based solvable variety. But the variety ${\cal B}_p{\cal B}_p\cap {\cal S}=({\cal B}_{p}\cap{\cal S}){\cal B}_p\cap\cal S$ is finitely based because by G. Higman's theorem (see\cite{HN}, 34.23), the
product of the nilpotent variety ${\cal B}_p\cap\cal S$ and the finitely based variety ${\cal B}_p$
must be finitely based.

Our statement on the basis rank also contrasts with the locally finite case. If nontrivial varieties
$\cal U$ and $\cal V$ consist of locally finite groups, then by Shmelkin's theorem \cite{S},
their product $\cal UV$ has finite basis
rank if and only if all  groups from $\cal U$ are nilpotent, all groups from $\cal V$ are abelian, and the
exponents of $\cal U$ and $\cal V$ are relatively prime.
\end{rk}

\section{Auxiliary finite  $p$-groups}\label{fin}

For an odd prime $p$ and  positive integers $\ell, m$, we will consider the sets $UT(\ell,p,m)$ of
$\ell\times \ell$ upper unitriangular matrices $A=(a_{ij})$, where for $j-i=t>0$,
the entry $a_{ij}$ is a $t$-linear function $V_{i}\times\dots \times V_{j-1}\to {\mathbb F}_p$
to the finite Galois field ${\mathbb F}_p$, where $V_1,\dots, V_{\ell-1}$ are $m$-dimensional
vector spaces over ${\mathbb F}_p$ (and $a_{ii}=1$).

Assume that $A=(a_{ij}), B=(b_{ij})\in UT(\ell,p,m)$. Then for $i\le j\le k$, the
component-wise product of the mappings $a_{ij}$ and $b_{jk}$ is a $(k-i)$-linear
mapping $V_{i}\times\dots V_{j-1}\times V_j\times\dots\times V_{k-1}\to {\mathbb F}_p$.
Therefore the matrix multiplication $AB=C$, where $c_{ik}=\sum_{j=1}^{\ell}a_{ij} b_{jk}$
is well defined in $UT(\ell,p,m)$. Moreover, it is associative for the ordinary reason,
and every matrix $A=I+B$ is invertible in $UT(\ell,p,m)$: $(I+B)^{-1}=I-B+...\pm B^{\ell-1}$ since $B$ is a
nil-triangle matrix. Hence the set
$UT(\ell,p,m)$ is a finite group under the matrix multiplication. Moreover, it is a $p$-group
since the number of possible entries $a_{ij}$ is a power of $p$ for every pair $(i,j)$ ($1\le i<j\le\ell$). Similarly, all
upper triangle matrices of this form (but the diagonal entries are any elements from ${\mathbb F}_p$
now) form an associative algebra over ${\mathbb F}_p$.

We will use the following properties of the groups $UT(p+1, p, m)$ (i.e. $\ell=p+1$
from now on) and of its subgroup $UT(p+1, p, m)^p$ generated by all the $p$th powers of  elements.

\begin{lemma} \label{U} (1) Let $A\in UT(p+1, p, m)$ and $C=A^p$. Then $c_{1,p+1}=a_{12}a_{23}\dots a_{p,p+1}$
and other entries of the matrix $C$ coincide with the corresponding entries of the identity matrix $I$.

(2) $UT(p+1, p, m)^p$ is a nontrivial central subgroup of exponent $p$.

(3) If $p\ge 3$ and $m\ge 4$, then $UT(p+1, p, m)^p$ has a nontrivial  cyclic subgroup containing no nontrivial
 element which can be presented in $UT(p+1, p, m)$
as a product of $m$ $p$th powers.
\end{lemma}
\proof (1) We have $A=I+B$, where $B$ is a nil-triangular matrix. Hence $A^p=\sum_{i=0}^p {p\choose i}B^i=I+B^p$
since ${p\choose i}$ is divisible by $p$ for $0<i<p$.  Thus $c_{1,p+1}=b_{12}b_{23}\dots b_{p,p+1}=a_{12}a_{23}\dots a_{p,p+1}$ and other entries of $B^p$ are $0$, which proves the statement.

(2) It follows from the claim (1) that there is a non-identity matrix $C$ in $UT(p+1, p, m)^p$.
It is straightforward that every matrix $C$ which is equal to $I$, except possibly in the entry $c_{1,p+1}$, belongs to
the center of $UT(p+1, p, m)$. Also we have $C^p=I$ for such a matrix since $C$ has zeros just
above the diagonal. This proves the statement (2).

(3) On one hand, each of the vector spaces $V_i$ has $p^m$ elements, and so there are $p^m$
linear functions on $V_i$ with values in ${\mathbb F}_p$. Hence the number of different
products $a_{12}\dots a_{p,p+1}$ does not exceed $(p^m)^p=p^{pm}$. It follows from the claim (1)
that the number of different products of the form $A_1^p\dots A_m^p$  in $UT(p+1, p, m)$
does not exceed $(p^{pm})^m=p^{pm^2}$.

On the other hand, the vector space of $p$-linear mappings $V_1\times\dots \times V_p\to {\mathbb F}_p$ has
dimension $m^p$ over ${\mathbb F}_p.$ This space is a linear envelope of $m^p$ linearly independent products of $p$
basic $1$-linear mappings from $V_i$ to ${\mathbb F}_p$ ($i=1,\dots,p$). Therefore the claim (1)
implies that the subgroup  $UT(p+1, p, m)^p$ has dimension $m^p$ over ${\mathbb F}_p,$
and so its order is $p^{m^p}$.

It follows that the subgroup $UT(p+1, p, m)^p$ has $(p^{m^p}-1)/(p-1)$ subgroups of order $p$.
Since the pairwise intersections of these cyclic subgroups are trivial, one of the subgroups
contains no nontrivial products of the form  $A_1^p\dots A_m^p$ if $(p^{m^p}-1)/(p-1)>p^{pm^2}.$
The last inequality follows from the inequalities $p\ge 3$ and $m\ge 4$, and the proof
is complete.

\endproof

\section{Construction of infinite torsion groups}\label{inf}

Our proof of Theorem \ref{main} makes use the approach from \cite{O89}. We recall a few definitions here.

Given a finite or infinite group alphabet ${\mathcal A}^{\pm 1}=\{a_1^{\pm 1}, a_2^{\pm 1},\dots\}$, we
we write $U\equiv V$ to express letter--by--letter equality of two words $U$ and $V$ over ${\mathcal A}^{\pm 1}$.
If $\mathcal A$ is a generating set of a group $G$, we
write $U=V$ whenever two words $U$ and $V$ over $\mathcal A ^{\pm
1}$ represent the same element of $G$; we identify the words over
$\mathcal A ^{\pm 1} $ and the elements of $G$ represented by them.

We shall assume that the series of groups $G(0), G(1), \dots$ is defined
according to the following scheme. By definition, $G(0)=F(\mathcal A)$ is
the free group with basis $\mathcal A$. We also set ${\cal R}_0={\cal X}_0=\emptyset.$

For $i\ge 1$, we assume by induction that the set of defining relators ${\cal R}_{i-1}$ of
the group $G(i-1)=\langle {\mathcal A}\mid {\cal R}_{i-1}\rangle$ is already defined,
as well as the sets ${\cal X}_j$ of {\it periods} of ranks $j<i$.

For $i\ge 1$, a word $X$ in the alphabet $\mathcal A ^{\pm 1} $ is called {\it
simple in rank} $i-1$, if it is not conjugate to a power of a
shorter word or to a power of a period of rank $\le i-1$ in the group $G(i-1)$. We
denote by $\mathcal X_i$ a maximal subset of words
satisfying the following conditions.

\begin{enumerate}
\item[1)] $\mathcal X_i$ consists of words of length $i$ which are simple
in rank $i-1$.

\item[2)] If $A, B\in \mathcal X_i$ and $A\not\equiv B$, then $A$ is not
conjugate to $B$ or $B^{-1}$ in the group $G(i-1)$.
\end{enumerate}

\noindent Each word from $\mathcal X_i$ is called a {\it period of
rank $i$.} For every period $A\in {\cal X}_i$ we fix an odd integer $n_A$
which should be large enough, i.e. $n_A\ge n_0$ for a large constant $n_0$.

Let $\mathcal S _i=\{ A^{n_A} \mid A\in \mathcal X _i\} $. Then the set $\mathcal R _i$ is defined by $\mathcal R _i=\mathcal R _{i-1} \cup \mathcal S _i$, and so $G(i)=\langle {\mathcal A}\mid {\cal R}_{i}\rangle$. By definition, the
group $G(\infty)$ is the inductive limit for the epimorphisms $G(0)\to G(1)\to\dots$, i.e.,
$G(\infty)=\langle {\mathcal A}\mid \cup_{i=0}^{\infty}{\cal R}_i\rangle.$

The group $G(\infty)$ depends on the choice of the exponents $n_A$, but there are some general   properties  for such constructions.

\begin{lemma}\label{peri} (1) (\cite{O89}, Theorem 26.4.3) Every element of $G(\infty)$ is  conjugate to
a power of a period of some rank $i\ge 1$.

(2) (\cite{O89}, Theorem 26.4.1) Every period $A$ has order $n_A$ in $G(\infty).$

\end{lemma}

\section{Main construction}\label{mc}

To combine the constructions from Sections \ref{fin} and \ref{inf} we fix a prime $p\ge n_0$, fix an integer $m\ge 4$ and
consider the finite $p$-group $U=UT(p+1,p,m)$. By Lemma \ref{U} (3), the subgroup $U^p$ has an element $g$ such
that no nontrivial element of $\langle g\rangle$ is a product of $m$ $p$th powers in $U$. Then we take a free
group $F$ admitting an epimorphism $f: F\to U$. Let $H=F^p$ be the subgroup generated by all $p$th powers in $F$.
Its image $f(H)=f(F^p)=f(F)^p= U^p$ is a central subgroup of $U$ by Lemma \ref{U} (2).

Being a free group itself, the subgroup $H$ has a free basis $\cal A$. So one can construct the groups
$G(0), G(1),\dots$ as in Section \ref{inf}, starting with $G(0)=H$. We will use the mapping $f$
to choose the exponents $n_A$ of the periods as follows.
If $A$ is a period of some rank $i\ge 1$ and
 $f(A)\in U^p\backslash \langle g\rangle$ or $f(A)=1$, then $n_A=p$;  otherwise $n_A=p^2.$

Let $N=N(m)$ be the kernel of the canonical mapping $H\to G(\infty)=\langle \cal A\mid \cal R\rangle$.

\begin{lemma} \label{ker} If two words $V$ and $W$ over ${\cal A}^{\pm 1}$ are equal in $G(\infty)$, then $f(V)=f(W)$.
\end{lemma}
\proof By the definition of the exponents $n_A$ and Lemma \ref{U} (2),  every relator $A^{n_A}$
belongs to $\Ker f$. Hence $N\le \Ker f$, and the statement follows.
\endproof

\begin{lemma} \label{ord}  Let $B$ be a word in the alphabet ${\cal A}^{\pm 1}$ and $B\ne 1$ in $G(\infty)$.

(1) If $f(B)\in U^p\backslash \langle g \rangle$
or $f(B)=1,$ then $B$ has order $p$ in $G(\infty)$.

(2) If $f(B)\in \langle g \rangle\backslash \{1\}$, then the order of $B$ in $G(\infty)$
is $p^2$.

\end{lemma}

\proof By Lemma \ref{peri} (1), $B$ is conjugate to a power of a period $A$ in $G(\infty)$;
say $B=ZA^kZ^{-1}$, and so the order of $B$ is equal to the order of $A^k$. The conditions
$f(B)\in U^p\backslash \langle g \rangle$
and $f(ZA^kZ^{-1})\in U^p\backslash \langle g \rangle$ are equivalent by Lemma \ref{ker}.
So are the conditions $f(B)=1$ and $f(ZA^kZ^{-1})=1$.

The condition $f(ZA^kZ^{-1})\in U^p\backslash \langle g \rangle$ is equivalent to $f(A^k)\in U^p\backslash \langle g \rangle $
by Lemma \ref{ker} . Similarly, the condition $f(B)\in \langle g \rangle\backslash \{1\}$
is equivalent to $f(A^k)\in \langle g \rangle\backslash \{1\}$.
Thus one may just assume that $B\equiv A^k.$

(1) Assume that $f(A^k)\in U^p\backslash \langle g \rangle$ or $f(B)=1$.
If $(p,k)=1$, then both $A^k$ and $A$ generate the
$p$-subgroup $\langle A \rangle.$ The same is true for the pair $f(A^k)$ and $f(A)$ since $U$ is a $p$-group.
Thus the order of $A^k$ is equal to the order of $A$ in $G(\infty)$, which is $p$
by Lemma \ref{peri} (2). If $p$ divides $k$, then the order of $A^k$ is $p$ since $A^{p^2}=1$ for any period $A$ and $A^k$ is nontrivial
in $G(\infty)$.

(2) Now suppose $f(A^k)\in \langle g \rangle\backslash \{1\}$. Since no element of $\langle g \rangle\backslash \{1\}$
is a $p$-th power (moreover, it cannot be a product  of $m$ $p$th powers in $U$), we have $(p,k)=1.$ Thus,
as above, one can replace $A^k$ by $A$ in the proof, and the order of $A$ is $n_A=p^2$ since
$f(A)\in \langle g \rangle\backslash \{1\}$.
 The lemma is proved.
\endproof

\begin{lemma} The subgroup $N=N(m)$ is normal in the group $F$.
\end{lemma}

\proof
 By definition, $N$ is the normal closure in $H$ of the powers $A^{n_A}$ of periods, where $n_A=p$ if $f(A)\in U^p\backslash \langle g\rangle$ or $f(A)=1$, and $n_A=p^2$ otherwise. By the statement (1), $N$ is the normal closure in $H$
of a bigger set $\cal Q$; namely,  $\cal Q$ consists of $p$-th powers of all the words $B$, where  $f(B)\in U^p\backslash \langle g \rangle$ or $f(B)=1,$ and of  $p^2$-th powers of all other words. Since $\langle g\rangle$ is a normal
subgroup of $U$ by Lemma \ref{U} (2), the set $\cal Q$ is invariant under conjugations in $F$. Therefore the
subgroup $N$ generated by $\cal Q$ is normal in $F$.
\endproof

Denote by $L(m)$ the quotient $F/N(m)$. Now we are prepared to prove the main lemma.

\begin{lemma} \label{FN} For given $m\ge 4$ the group $L(m)$ does not belong to the variety ${\cal B}_p{\cal B}_p$ but it satisfies the identical relation $(x_1^p\dots x^p_m)^p=1.$
\end{lemma}
\proof The group $H=F^p$ has an element $B$ such that
$f(B) \in \langle g\rangle\backslash\{1\}.$ Its order in $H/N\le F/N=L(m)$ is
$p^2$ by Lemma \ref{ord} (2). Therefore the verbal subgroup $(F/N)^p$ of
$F/N$ does not belong to the variety ${\cal B}_p$, and so
$F/N\notin {\cal B}_p{\cal B}_p$.

If $h=g_1^p\dots g_m^p$ for some $g_1,\dots, g_m\in F $,
then $f(h)=f(g_1)^p\dots f(g_m)^p\in U^p.$
So either $f(h)=1$ or $f(h)\notin \langle g \rangle$ by the definition
of the element $g$. Hence $h^p=1$ in $H/N$ by Lemma \ref{ord} (1).
Since $(g_1^p\dots g_m^p)^p=1$ in $F/N$ for every $m$-tuple $(g_1,\dots, g_m)$,
the required identity is obtained.
\endproof

{\bf Proof of Theorem \ref{main}.} If the system (\ref{BpBp}) is equivalent to
a finite system of identities, then it is equivalent to its own finite subsystem

\begin{equation}\label{m}
x_1^{p^2}=1, (x_1^px_2^p)^p=1,\dots, (x_1^px_2^p\dots x_m^p)^p=1
\end{equation}
for some $m$, and one may assume that $m\ge 4$. Note that the last identity implies
all the preceding ones. By Lemma \ref{FN}, there is a group $L(m)$ satisfying
the system (\ref{m}) but not satisfying the system (\ref{BpBp}). Thus, the theorem is
proved by contradiction.
$\Box$

\bigskip

{\bf Proof of part (1) in Theorem \ref{rank}.} It suffices for every $k\ge 1$, to find a group $M(k)$
which does not belong to the product ${\cal B}_p{\cal B}_p$ but every $k$-generated subgroup
of $M(k)$ lies in ${\cal B}_p{\cal B}_p$.

Note that for given $k$, all $k$-generated subgroups $K$ of the groups $UT(p+1, p, m)$ ($m=1,2,\dots$)
have bounded orders since their nilpotency classes, exponents and the numbers of generators are bounded
(by $p$, $p^2$ and $k$, respectively). Therefore there exists $m=m(k)$ such that every element of $K^p$
is a product of at most $m$  $p$-th powers in every such a subgroup $K$.

Let now $M(k)=G(\infty)$ be the group whose definition in Section  \ref{mc} depends on the group $U=U(p+1,p,m)$
and the choice of $g\in U^p$, where nontrivial elements of $\langle g \rangle$ are not products of $ m$
$p$-th powers in $U$. It follows that $K^p\cap \langle g\rangle =\{1\}$ for every $k$-generated subgroup $K$
of $U$.

Assume now that $L$ is a $k$-generated subgroup of $M(k)$. Since $f(L)$ is a $k$-generated subgroup of $U$,
 we have $f(L)^p\cap \langle g\rangle =\{1\}$. Therefore by Lemma \ref{ord}, we have $z^p=1$ for every element $z$ of the subgroup $L^p$, that is, the subgroup $L$ satisfies the system of identities (\ref{BpBp}). Since by Lemma \ref{FN}, the group $M(k)$ itself does not belong to the variety  ${\cal B}_p{\cal B}_p$, the proof of the statement (1) is complete.

\section{Products of varieties and their basis rank}

For the second part of Theorem \ref{rank}, we will use the concept introduced in \cite{BNNN}. A group $G$ is
said to be {\it discriminated} by a group $D$ if for every finite subset $S\subset G$, there is a homomorphism
$G\to D$ injective on $S$. If a free group of infinite countable rank $F_{\infty}(\cal V)$ of a variety $\cal V$
is discriminated by a group $D\in \cal V$, then $var D=\cal V$, i.e. the variety $\cal V$ is generated by
the group $D$. Clearly a group $C$ discriminates a group $G$ if a subgroup $D\le C$ discriminates $G$.
If a factor group $B$ of $D$ discriminates the group $F_{\infty}(\cal V)$ and $D\in \cal V$, then $D$
discriminates the group $F_{\infty}(\cal V)$ too since every homomorphism $F_{\infty}({\cal V})\to B$
lifts to a homomorphism $F_{\infty}({\cal V})\to D$ by the universal property of $\cal V$-free groups.

\begin{theorem} \label{UV} Let a group $G$ be $m$-generated, $m\ge 2$, such that $G$ belongs to a variety $\cal V$.
If $G$ contains a finitely generated subgroup $D$ of infinite index which discriminates the $\cal V$-free
group $F_{\infty}(\cal V)$, then for an arbitrary variety of groups $\cal U$, the $\cal UV$-free group
$F_{\infty}(\cal UV)$ is discriminated by an $m$-generated group from $\cal UV.$
Therefore the basis rank of the product $\cal UV$ does not exceed $m$.
\end{theorem}

\proof One may assume that $\cal V \ne \cal O$, where $\cal O$ is the variety of all groups, because $\cal UO=O$ and
the (absolutely) free group $F_m$ of rank $m$ discriminates the group $F_{\infty}=F_{\infty}(\cal O)$ since $F_m$
contains an isomorphic copy of $F_{\infty}$.

Let $N$ be a normal subgroup of $F_m$ such that $F_m/N\cong G$. We have $N\ne \{1\}$ since $G\in \cal V\ne O.$
The group $F_m/U(N)$, where $U(N)$ is the verbal subgroup of $N$ corresponding to the variety $\cal U$, obviously
belongs to the variety $\cal UV$. It follows from the assumptions that there is a subgroup $K\le F_m$
such that $N\le K$, $K/N\cong D$ and the index $[F_m: K]$ is infinite.

The group $K$ is free, being a subgroup of $F_m$, and $K$ is not finitely
generated since it has infinite index in $F_m$ and contains a nontrivial normal subgroup $N$ of $F_m$ (see \cite{LS}, Proposition I.3.11).
Thus we have an infinite free basis $y_1, y_2,\dots $ in $K$.

Since $K/N$ is a finitely generated group, there is a finite subset $y_1,\dots, y_s$ generating the subgroup
$K$ modulo $N$. Therefore there are words $v_{s+1}=v_{s+1}(y_1,\dots, y_s), v_{s+2}=v_{s+2}(y_1,\dots, y_s),\dots$
such that the products $z_1=y_{s+1}v_{s+1},z_2=y_{s+2}v_{s+2},\dots$ belong to $N$. The set $y_1,\dots, y_s, z_1,  \dots$ is also a
free basis of $K$, since it is obtained from a basis by Nielsen transformations.

Now we consider the auxiliary wreath product $W=F_{\infty}({\cal U}) wr D$ which belongs to the product $\cal UV$ and discriminates the group $F_{\infty}(\cal UV)$
(see \cite{Ba} or \cite{HN}, Corollary 22.44). Recall that $W$
is generated by a copy $A$ of $F_{\infty}(\cal U)$ and by $D$, $W$ is a semidirect product of the normal closure
$M$ of $A$ in $W$ and the subgroup $D$, $M$ is isomorphic to a direct power
of $A$, and so $M\in \cal U$.

Let $\varepsilon$ be the canonical homomorphism $K\to K/N=D.$ To construct  an epimorphism $\eta: K\to W$, we will
modify $\varepsilon$ as following. The elements of the free basis $y_1,...,y_s$ of $K$ are mapped by $\eta$ into
$D$ according to $\varepsilon$ but the set $\{z_1, z_2,\dots\}$ is mapped onto a generating set of $A$.
The canonical projection $\alpha$ of $W$ onto the semidirect factor $D$ removes the difference between $\eta$
and $\varepsilon$, that is, $\alpha\eta(x)=\varepsilon(x)$ for every $x\in \{y_1,\dots,y_s, z_1,\dots\}$.
It follows that $\eta(N)\le \ker \alpha=M$.

Since $\eta(N)\le M\in \cal U$, we have $\eta(U(N))=U(\eta(N))=\{1\}$. Hence the epimorphism $\eta$ factors
through an epimorphism $K/U(N)\to W$. Since $W$ discriminates the group $F_{\infty}(\cal UV)$, so do
the group $K/U(N)$ and the $m$-generated group $F_m/U(N)$ containing $K/U(N)$. This completes the proof of Theorem \ref{UV}.
\endproof

\begin{cy} \label{UB} For any variety $\cal U$ and sufficiently large odd integer $n$, the variety ${\cal UB}_n$ has
basis rank $2$.
\end{cy}

\proof
Recall that for large enough odd $n$ the group
$F_{\infty}({\cal B}_n)$ is isomorphic to a subgroup of $F_2({\cal B}_n)$ (see \cite{Sh} or \cite{O89},
Corollary 35.6). Hence $F_{\infty}({\cal B}_n)$ is discriminated by the $2$-generated group
$F_2({\cal B}_n)$, and the group $F_2({\cal B}_n)$ contains its own isomorphic copy as a subgroup
of infinite index. Therefore the basis rank of the product ${\cal U}{\cal B}_n$ is at most $2$
by Theorem \ref{UV}. Clearly, it is greater than $1$, since this product is not an abelian variety.
\endproof

{\bf Proof of part (2) in Theorem \ref{rank}}. The statement follows from Corollary \ref{UB} with ${\cal U=\cal B}_p$ and $n=p$. $\Box$.

\begin{rk} One can replace ``large enough odd'' by ''divisible by a large enough power of $2$'' in the formulation
of Corollary \ref{UB} since under the latter assumption, the group $F_{\infty}({\cal B}_n)$ also embeds into
the group $F_2({\cal B}_n)$ \cite{IO}.
\end{rk}

\begin{cy} \label{m1} If a variety $\cal V$ is nontrivial and the group $F_{\infty}(\cal V)$ is discriminated
by an $m$-generated group $D$ from $\cal V$, then for every variety $\cal U$, the group $F_{\infty}(\cal UV)$
is discriminated by an $(m+1)$-generated group from $\cal UV$. In particular, the basis rank of the
product $\cal UV$ does not exceed $m+1$.
\end{cy}

\proof The group $D$ is a homomorphic image of the $\cal V$-free group $F_m(\cal V)$, and so $F_m(\cal V)$
discriminates the group $F_{\infty}(\cal V)$ as well. It is naturally embedded in $F_{m+1}(\cal V)$ as
a subgroup of infinite index. Indeed if we had $[F_{m+1}({\cal V}):F_m (\cal V)]<\infty$ and $(a_1,\dots a_{m+1})$
was a $\cal V$-free basis of $F_{m+1}(\cal V)$, then the endomorphism $a_1\mapsto 1, a_2\mapsto a_1,\dots a_{m+1}\mapsto a_m$
would show that $[F_{m}({\cal V}):F_{m-1}(\cal V)]<\infty$ as well.  By induction, it would follow that $|F_{m+1}({\cal V})|<\infty$, contrary to
the fact that a nontrivial discriminating group cannot be finite (\cite{HN}, 17.32). Hence the corollary
follows from Theorem \ref{UV} applied to $F_{m+1}({\cal V})$.
\endproof

\begin{rk} On one hand, for any two group varieties $\cal U$ and $\cal V$ with $\cal U\ne O$, where $\cal O$ is the
variety of all groups, the basis rank of $\cal V$ does not exceed the basis rank of the product $\cal UV$ (\cite{HN},
25.12).  On the other hand,  the free abelian group of any rank is discriminated by an infinite cyclic group (\cite{HN}, 17.6). Hence $m=1$
if $\cal V$  is the variety of all abelian groups. But for any nontrivial variety $\cal U$, the
product $\cal UV$ is not abelian, and so its basis rank is at least $2$. Therefore the statement of
Corollary \ref{m1} fails when replacing $m+1$ by $m$ in its formulation.

Note that the weaker formulation
with $2m+1$ instead of $m+1$ can be found in \cite{BNNN} and \cite{HN}, Corollary 25.22. H.Neumann mentioned
in \cite{HN} that she needed additional assumptions to get $(m+1)$-generated group discriminating the group
$F_{\infty}(\cal UV)$. Since Corollary \ref{m1} gives the estimate $m+1$ without any additional restrictions,
Theorem \ref{UV} is of certain interest independently of the application to Theorem \ref{rank}.
\end{rk}

\bigskip

{\bf Nicholas S. Boatman: } Vanderbilt University, 2012. {\it E-mail}: nickboatman@yahoo.com

\medskip

{\bf Alexander A. Olshanskii:}  Vanderbilt University, Nashville
37240, USA, and Moscow State University, Moscow 119991, Russia.
{\it E-mail}: alexander.olshanskiy@vanderbilt.edu

\end{document}